\documentclass[11pt]{article}
\usepackage{latexsym,amssymb}
\usepackage{amsmath}
\usepackage[mathscr]{eucal}
\usepackage[abbrev]{amsrefs}

\newtheorem{Theorem}{\bf Theorem}[section]
\newtheorem{Lemma}{\bf Lemma}[section]
\newtheorem{Proposition}{\bf Proposition}[section]  
\newtheorem{Corollary}{\bf Corollary}[section]
\newtheorem{Remark}{\bf Remark}[section]
\newtheorem{Example}{\bf Example}[section]
\newtheorem{Definition}{\bf Definition}[section]
\newtheorem{Condition}{\bf Condition}[section]

\newenvironment{theorem}{\begin{Theorem}$\!\!\!$}{\end{Theorem}}
\newenvironment{lemma}{\begin{Lemma}$\!\!\!$}{\end{Lemma}}
\newenvironment{proposition}{\begin{Proposition}$\!\!\!$}{\end{Proposition}}

\newenvironment{remark}{\begin{Remark}$\!\!\!$}{\end{Remark}}
\newenvironment{definition}{\begin{Definition}$\!\!\!$}{\end{Definition}}

\numberwithin{equation}{section}
\begin{document}
\title{New characterizations of log-concavity}
\author{Kazuhiro Ishige, Paolo Salani and Asuka Takatsu}

\date{}
\maketitle
\begin{abstract}
We introduce a notion of $F$-concavity which largely generalizes the usual concavity. 
By the use of the notions of closedness under positive scalar multiplication and closedness under positive exponentiation 
we characterize power concavity and power log-concavity 
among nontrivial $F$-concavities,  respectively. In particular, we have a characterization of log-concavity as the only $F$-concavity which is closed both under positive scalar multiplication and positive exponentiation.
Furthermore, we discuss the strongest $F$-concavity 
preserved by the Dirichlet heat flow, characterizing log-concavity also in this connection.
\end{abstract}
\vspace{40pt}
\noindent Addresses:

\smallskip
\noindent 
K. I.: Graduate School of Mathematical Sciences, The University of Tokyo, 3-8-1 
Komaba, Meguro-ku, Tokyo 153-8914, Japan\\
\noindent 
E-mail: {\tt ishige@ms.u-tokyo.ac.jp}\\

\smallskip
\noindent 
P. S.: DiMaI Dipartimento di Matematica e Informatica``U. Dini", 
Universit\`a di Firenze, viale Morgagni 67/A, 50134 Firenze\\
\noindent 
E-mail: {\tt paolo.salani@unifi.it}\\

\smallskip
\noindent 
A. T.: Department of Mathematical Sciences,  Tokyo Metropolitan University, 1-1
Minami-osawa, Hachioji-shi, Tokyo 192-0397, Japan\\
\noindent 
E-mail: {\tt asuka@tmu.ac.jp}\\

\vspace{20pt}
\vspace{3pt}
\newpage
%%%%%%%%%%%%%%%%%%%%%%%%%%%%%%%%%%%%
\section{Introduction}
%%%%%%%%%%%%%%%%%%%%%%%%%%%%%%%%%%%%
\subsection{Prologue}
Concavity and its variations are useful notions to describe the shape of functions 
and they have fascinated many mathematicians. In particular, in the theory of partial differential equations, an interesting issue is to determine if the solution to a Dirichlet problem in a convex domain shares some concavity property.
Pioneering works in this direction are the following results.
\begin{itemize}
  \item[(A)] Gabriel \cite{G}\\
  Every superlevel set of the Newtonian potential of a bounded convex set in ${\bf R}^3$ is convex.
  \item[(B)] Makar-Limanov \cite{ML}\\
  The square root of the torsion function of a planar bounded convex domain is concave.
  \item[(C)] Brascamp and Lieb~\cite{BL}\\[-20pt]
  \begin{itemize}
  \item[(i)]
  The Dirichlet heat flow in a convex domain preserves logarithmic concavity. 
  \item[(ii)]
  The first positive Dirichlet eigenfunction of the Laplacian in a bounded convex domain is logarithmically concave.
  \end{itemize}
\end{itemize}
All these results have been generalized in several ways, 
and they can be conveniently expressed 
in terms of {\em power concavity} (see Definition \ref{Example:1.1} below).
In particular, assertions~(A), (B) and (C) can be rephrased as follows.
\begin{itemize}
  \item[(A')] 
  The Newtonian potential of a bounded convex set is {\em quasiconcave}.
  \item[(B')] 
  The torsion function of a bounded convex domain is {\em $1/2$-concave}.
  \item[(C')] 
  \begin{itemize}
  \item[(i)]
  The Dirichlet heat flow in a convex domain preserves {\em  $0$-concavity}. 
  \item[(ii)]
  The first positive Dirichlet eigenfunction of the Laplacian in a bounded convex domain is {\em $0$-concave}.
  \end{itemize}
\end{itemize} 
A quasiconcave (or $(-\infty)$-concave) function is a function whose superlevel sets are all convex (see later for more details), while the definition of $p$-concavity for $p\in{\bf R}$ is as follows.
\begin{Definition}
\label{Example:1.1} {\rm (Power concavity)} 
Let $p\in {\bf R} $ and define
$$
\Phi_p(\tau):=
\left\{
\begin{array}{ll}
\displaystyle{\frac{\tau^p-1}{p}} & \quad\mbox{for}\quad p\not=0,\vspace{5pt}\\
\log\tau& \quad\mbox{for}\quad p=0,
\end{array}
\right.
$$
for $\tau>0$. Set $\Phi_p(0):=\lim_{\tau\to +0}\Phi_p(\tau)$. 
Let $\Omega$ be a convex domain in ${\bf R}^N$.
A nonnegative function $f$ in $\Omega$ is said {\em $p$-concave} in $\Omega$ if 
$$
 \Phi_p(f(1-\mu)x+\mu y))\ge (1-\mu)\Phi_p(f(x))+\mu \Phi_p(f(y))
$$
for $x$, $y\in\Omega$ and $\mu\in[0,1]$. 
\end{Definition}
Roughly speaking, 
$f$ is $p$-concave in $\Omega$ if:
\begin{center}
$f^p$ is concave, for $p>0$;\quad
$\log f$ is concave, for $p=0$;\quad
$f^p$ is convex, for $p<0$.
\end{center}
For $p=1$, we clearly get back to the usual concavity.
A $0$-concave function is commonly called {\em log-concave}.
Log-concavity deserves a special place among power concavities, due to its interplay with the Dirichlet heat flow 
(as we will see also in this paper) and its relevance in many fields of mathematics  and in other applied sciences (see for instance \cite{An, Co, SW} for presentations of log-concave functions from different perspectives).

Power concavity is the most common variation of the usual concavity. In this paper we further generalize this notion and we introduce a notion of $F$-concavity, embracing power concavity and other interesting variations of concavity. Then we
characterize different concavities (in particular, log-concavity) through some of their relevant properties and in connection 
with the Dirichlet heat flow. 

\subsection{Notations and definitions}
Unless otherwise stated, we denote by $I$ an interval in ${\bf R}$ with $\mbox{Int}\,I\not=\emptyset$, 
where $\mbox{Int}\,I$ is the interior of $I$. Let $\Omega$ be a convex domain in ${\bf R}^N$. We set
\begin{equation*}
\begin{split}
{\mathcal A}_\Omega(I):=
 & \{f\,:\,\mbox{$f$ is a function in $\Omega$ such that $f(x)\in I$ for $x\in\Omega$}\},\\
{\mathcal A}^*_\Omega(I):=
 & \{f\,:\,\mbox{$f$ is a function in $\Omega$ such that $\inf_\Omega f\in \mbox{Int}\,I$ and $\sup_\Omega f\in \mbox{Int}\,I$}\}\,,\\
 \mathcal{C}_\Omega^+:= &\{f\,:\,\mbox{$f$ is a constant function in $\Omega$}\}.
\end{split}
\end{equation*}
Note that ${\mathcal A}^*_\Omega(I)\subset{\mathcal A}_\Omega(I)$. 

\begin{definition}
\label{Definition:1.0} Let $I$ be an interval in ${\bf R}$. A function $F:I\to{\bf R}\cup\{-\infty\}$ is said {\em admissible} on $I$ if 
$F\in C(\mbox{Int}\,I)$ and $F$ is strictly increasing on $I$. 
\end{definition}
\begin{definition}
\label{Definition:1.1}
Let $\Omega$ be a convex domain in ${\bf R}^N$ and $I$ an interval in ${\bf R}$.
\begin{itemize}
  \item[{\rm (a)}] 
  Let $F$ be admissible on $I$. For any $f\in \mathcal{A}_\Omega(I)$, 
  we say that $f$ is {\em $F$-concave} in $\Omega$ if
  \begin{equation}\label{eq:F-concavity}
  F(f((1-\mu)x+\mu y))\ge (1-\mu)F(f(x))+\mu F(f(y))
  \end{equation}
  for $x$, $y\in\Omega$ and $\mu\in[0,1]$. 
  We denote by $\mathcal{C}_\Omega[F]$ the set of $F$-concave functions in $\Omega$. 
  \item[{\rm (b)}] 
 Let $F_1$ and $F_2$ be admissible on $I$. 
 We say that {\em $F_1$-concavity is stronger ({\it resp.} weaker) than $F_2$-concavity} in $\Omega$ if 
 $ \mathcal{C}_\Omega[F_1]\subseteq\mathcal{C}_\Omega[F_2]
  \,\,
  (\text{resp. }\mathcal{C}_\Omega[F_2]\subseteq\mathcal{C}_\Omega[F_1]). $
   \item[{\rm (c)}] We say that $F$-concavity is {\em trivial} in $\Omega$ if 
   $\mathcal{C}_\Omega[F]=\mathcal{A}_\Omega(I)\cap\mathcal{C}_\Omega^+$.
\end{itemize}
\end{definition}
Notice that, with the definition above, 
we can compare only concavities admissible on the same interval
(see Section~4.2).
\begin{remark}\label{rem:mean}
For an admissible function $F$, we can equivalently rewrite \eqref{eq:F-concavity} as follows:
$$
f((1-\mu)x+\mu y)\ge F^{-1}\left((1-\mu)F(f(x))+\mu F(f(y))\right).
$$
The right-hand side of the above inequality is called the {\em $F$-mean} of $f(x)$ and $f(y)$.
For properties of $F$-means, we refer to e.g. \cite{HLP, Na}. 
\end{remark}
\begin{remark}
\label{rem-F}
Let $F$ be admissible on an interval $I$. 
\begin{itemize}
  \item[{\rm (i)}] 
  Set $\tilde{F}(\tau):=-F(-\tau)$ for $-\tau\in I$. 
  Then $\tilde{F}$ is admissible on $-I$ and 
  $\mathcal{C}_\Omega[\tilde{F]}=-\mathcal{C}_\Omega[F]$. 
  \item[{\rm (ii)}] 
  For any $\lambda>0$, set $\tilde{F}(\tau)=F(\lambda\tau)$ for $\tau\in\lambda^{-1}I$. 
  Then $\tilde{F}$ is admissible on $\lambda^{-1}I$ and 
  $\mathcal{C}_\Omega[\tilde{F}]=\lambda^{-1}\mathcal{C}_\Omega[F]$. 
  \item[{\rm (iii)}]
   If $F$ is admissible on an interval $I$,  then $\tilde{F}(\tau)=F(\tau)+c$ is admissible on  $I$ 
    for every $c\in{\bf R}$ and $\mathcal{C}_\Omega[\tilde{F}]=\mathcal{C}_\Omega[F]$.
  \item[{\rm (iv)}] 
  Let $F$ be admissible on an interval $I$.
  Then the restriction $\tilde{F}$ of $F$ to any subinterval $J\subseteq I$ is admissible on $J$ and 
  $\mathcal{C}_\Omega[\tilde{F}]=\mathcal{C}_\Omega[F]\cap\mathcal{A}_\Omega(J)$.
  \end{itemize}
\end{remark}
Observe that, 
for admissible functions~$F_1$ and $F_2$ on an interval $I$,  
if ${\mathcal C}_\Omega[F_1]={\mathcal C}_\Omega[F_2]$, then 
\begin{equation}
\label{eq:1.1}
\qquad
{\mathcal C}_{\lambda P\Omega+z}[F_1]={\mathcal C}_{\lambda P\Omega+z}[F_2]
\quad\mbox{for}\quad\lambda>0,\,\,\, P\in O(N),\,\,\,z\in{\bf R}^N.
\end{equation}
Here $O(N)$ is the orthogonal group in dimension~$N$. 

Before giving some examples of interesting $F$-concavities, 
we introduce the following notions, which are crucial to the scope of this paper.
\begin{definition}
\label{Definition:closeness}
Let $\Omega$ be a convex domain in ${\bf R}^N$ and $F$ admissible on an interval $I$. 
\begin{itemize}
 \item[{\rm (a)}] 
 We say that $F$-concavity is {\em closed under positive scalar multiplication} 
 if for any $f\in\mathcal{C}_\Omega[F]$, it holds $\lambda f\in\mathcal{C}_\Omega[F]$ for every $\lambda>0$.
 \item[{\rm (b)}] 
 We say that $F$-concavity is {\em closed under positive exponentiation} if
 for any $f\in\mathcal{C}_\Omega[F]$, it holds $f^r\in\mathcal{C}_\Omega[F]$ for every $r>0$.
\end{itemize}
\end{definition}
Clearly, properties (a) and/or (b) impose some restrictions about the interval $I$; we leave to the reader the discussion of these restrictions.

\begin{Example}
\label{Example:1.0}
The function $\Phi_p$ is admissible on $I=[0,\infty)$ for $p\in{\bf R}$ and power concavities are of course a relevant example of $F$-concavities.
Let us recall some properties of $p$-concave functions:
\begin{itemize}
  \item[{\rm (i)}] 
  If $p<q$, then $q$-concavity is stronger than $p$-concavity.
  \item[{\rm (ii)}] 
  $p$-concavity is closed under positive scalar multiplication, that is, 
  if $f$ is $p$-concave, then $\lambda f$ is also $p$-concave for $\lambda>0$.
  \item[{\rm (iii)}] 
  If $f$ is $p$-concave and $r>0$, then $f^r$ is $p/r$-concave. 
  In particular, log-concavity is closed under positive exponentiation, that is, 
  if $f$ is log-concave, then $f^r$ is also log-concave for $r>0$.
 \end{itemize}
 \end{Example}

\begin{Example}
\label{Example:1.2}
Let $p\in {\bf R}$. 
Set $\Phi_p^*(0):=-\infty$ if $p>0$ and otherwise $\Phi_p^*(\tau):=\Phi_p(\tau)$. 
Then $\Phi_p^*$ is admissible on $I=[0,\infty)$ such that $\Phi_p^*(0)=-\infty$. 
It easily follows that 
$$
{\mathcal C}_\Omega[\Phi_p]={\mathcal C}_\Omega[\Phi_p^*]\quad\mbox{if}\quad p\le 0,
\qquad
{\mathcal C}_\Omega[\Phi_p]\subsetneq{\mathcal C}_\Omega[\Phi_p^*]\quad\mbox{if}\quad p>0.
$$
Indeed, if $p>0$, then  
$$
\chi_K\in{\mathcal C}_\Omega[\Phi_p^*]\setminus {\mathcal C}_\Omega[\Phi_p], 
$$
where $\chi_K$ is the characteristic function of a convex set $K\subsetneq\Omega$. 
$\Phi_p^*$-concavity has been used in the study of the heat flow in ${\bf R}^N$ {\rm ({\it see e.g.} \cite{BL})}. 
\end{Example}
\begin{Example}
\label{Example:1.3} {\rm (Power log-concavity)} 
Let $\alpha \in {\bf R}$ and define
$$
L_\alpha(\tau):=
-\Phi_\alpha(-\log\tau)=
\left\{
\begin{array}{ll}
-\displaystyle{\frac{1}{\alpha}}\left[(-\log \tau)^\alpha-1\right]\quad & \mbox{if}\quad\alpha\not=0,\vspace{7pt}\\
-\log(-\log\tau)\quad & \mbox{if}\quad\alpha=0,
\end{array}
\right.
$$
for $\tau\in(0,1)$. Set $L_\alpha(0):=\lim_{\tau\to +0}L_\alpha(\tau)$ 
and $L_\alpha(1):=\lim_{\tau\to 1-0}L_\alpha(\tau)$. Then $L_\alpha$ is admissible on $I$, where $I=[0,1]$ if $\alpha > 0$ and $I=[0,1)$ if $\alpha \leq 0$.
Notice that $L_\alpha(0)=-\infty$ if and only if $\alpha\ge 0$ $($see also Theorem~{\rm\ref{Theorem:1.1}~(b)}$)$.
For $f \in \mathcal{A}_\Omega(I)$,
we say that $f$ is {\em $\alpha$-log-concave} in $\Omega$ if $f$ is $L_\alpha$-concave in $\Omega$.
The following properties hold
{\rm ({\it see} \cite[{\it Section}~2]{IST})}.
\begin{itemize}
 \item[{\rm (i)}] 
 Let $0\le f(x)\le 1$ in $\Omega$. Then $f$ is log-concave in $\Omega$ 
 if and only if $f$ is $1$-log-concave in~$\Omega$.
 \item[{\rm (ii)}] 
 If $\alpha<\beta$, then $\alpha$-log-concavity is stronger than $\beta$-log-concavity.
 \item[{\rm (iii)}] 
 If $\alpha\le 1$ and $f$ is $\alpha$-log-concave, 
 then $\lambda f$ is also $\alpha$-log-concave for $0<\lambda < (\sup_\Omega f)^{-1}$.
 \item[{\rm (iv)}]
 $\alpha$-log-concavity is closed under positive exponentiation 
 {\rm({\it see also Example}~\ref{Example:1.1}~(ii))}.
 \item[{\rm (v)}]
 For $1/2\le\alpha\le 1$, $\alpha$-log-concavity is preserved by the Dirichlet heat flow in $\Omega$ 
 {\rm ({\it see also Proposition}~\ref{Proposition:4.2})}.
\end{itemize}

Set
$$
f_\alpha(x):=
\left\{
\begin{array}{ll}
\exp\left(-\left(1+\alpha|x|\right)^{1/\alpha}\right) & \mbox{if}\quad\alpha\not=0,\vspace{5pt}\\
\exp(-\exp(|x|)) & \mbox{if}\quad\alpha=0,
\end{array}
\right.
$$
for $x\in{\bf R}^N$  if $\alpha\geq 0$ and otherwise  for $x\in B(0,1/|\alpha|)$.
Here $B(0,r):=\{x\in{\bf R}^N\,:\,|x|<r\}$ for $r>0$. 
Then 
$$
L_\alpha(f_\alpha(x))=-|x|\quad\text{for every }\alpha\in{\bf R}
$$
and $f_\alpha$ is $\alpha$-log-concave in ${\bf R}^N$ if $\alpha\geq 0$ and in $B(0,1/|\alpha|)$ otherwise.

Notice also that the function $\exp(-|x|^2)$ is $\alpha$-log-concave in ${\bf R}^N$ if and only if $\alpha\geq1/2$. Similarly the function $a\exp(-b|x|^2)$ is $(1/2)$-log-concave in ${\bf R}^N$ for  $0\leq a\leq1$ and $b\geq 0$.
\end{Example}

As already said, a function whose superlevel sets are all convex is said 
{\em quasiconcave} or {\em $(-\infty)$-concave}. 
More explicitly, $f$ is said quasiconcave in $\Omega$ if 
$$
f((1-\mu)x+\mu y)\ge \min\{f(x),\,f(y)\}
$$
for $x$, $y\in\Omega$ and $\mu\in[0,1]$. 
We denote by $\mathcal{C}_\Omega^-$ the sets of quasiconcave functions in $\Omega$. 
It easily follows that
\begin{equation}
\label{eq:1.3}
\mathcal{A}_\Omega(I)\cap\mathcal{C}_\Omega^+
\subseteq\mathcal{C}_\Omega[F]
\subseteq\mathcal{C}_\Omega^-
\end{equation} 
for every admissible functions $F$ on $I$. 
The second inclusion in \eqref{eq:1.3}, roughly speaking, 
tells that quasiconcavity is the weakest among conceivable concavities. 
\begin{remark}
\label{Remark:1.2}
Quasiconcavity may be regarded as the limit of $p$-concavity as $p\to-\infty$ and 
for this reason it is also denoted as $(-\infty)$-concavity.
 Indeed, for $p\neq 0$, the $p$-concavity of a positive function $f$ can be rewritten as follows 
 {\rm ({\it see also Remark}~\ref{rem:mean})}:
$$
f((1-\mu)x+\mu y)\ge \left[(1-\mu)f(x)^p+\mu f(y)^p\right]^{1/p}:=M_p(f(x),f(y);\mu)\,
$$
and $\lim_{p\to-\infty}M_p(a,b;\mu)=\min\{a,b\}$ 
for $a$, $b>0$ and $\mu\in(0,1)$. 

Despite this, quasiconcavity can not be expressed as $F$-concavity for any admissible $F$ 
and the following inclusion is indeed strict:
$$
\bigcup_{\mbox{admissible F}}\mathcal{C}_\Omega[F]\subsetneq\mathcal{C}_\Omega^{-}. 
$$
For instance, the function $f$ in ${\bf R}$ defined by 
$$
f(x):=\left\{
\begin{array}{ll}
1 & \mbox{for}\quad x\in(0,1],\\
2 & \mbox{for}\quad x\in(1,2),\\
0 & \mbox{otherwise},
\end{array}
\right.
$$
is quasiconcave in ${\bf R}$ but is not $F$-concave in ${\bf R}$ for any admissible $F$ on $[0,2]$. 

Notice also that quasiconcavity is closed under positive scalar multiplication. 
\end{remark}

\subsection{Main results}
We state hereafter the main results of this paper. 
First we specify trivial $F$-concavities, 
then we identify nontrivial $F$-concavities closed under positive scalar multiplication and positive exponentiation, respectively.
Finally, we characterize log-concavity in connection with the Dirichlet heat flow.
\begin{theorem}
\label{Theorem:1.1}
Let $F$ be admissible on an interval $I$. 
\begin{itemize}
  \item[{\rm (a)}] 
  Let $\Omega$ be a convex domain in ${\bf R}^N$ with $\Omega\not={\bf R}^N$. 
  Then $F$-concavity is nontrivial in $\Omega$. 
  \item[{\rm (b)}] 
  $F$-concavity is trivial in ${\bf R}^N$ if and only if $\,\inf_I F>-\infty$. 
  Furthermore, 
  if $F$ is nontrivial and $f\in{\mathcal C}_{{\bf R}^N}[F]\setminus {\mathcal C}^+_{{\bf R}^N}$, then $\inf_{{\bf R}^N} F(f)=-\infty$.
\end{itemize} 
\end{theorem}
It follows from Theorem~\ref{Theorem:1.1}~(b) that 
$\Phi_p$-concavity with $p>0$ and $L_\alpha$-concavity with $\alpha<0$ are trivial in ${\bf R}^N$.
\begin{theorem}
\label{Theorem:1.2}
Let $\Omega$ be a convex domain in ${\bf R}^N$.
Let $F_1$ and $F_2$ be admissible on an open interval $I$ 
such that $F_1$-concavity is nontrivial in $\Omega$. 
Then ${\mathcal C}_\Omega[F_1]={\mathcal C}_\Omega[F_2]$ if and only if 
there exists a pair $(A,B)\in (0,\infty)\times{\bf R}$ such that 
\begin{equation}
\label{eq:1.4}
F_1(\tau)=AF_2(\tau)+B\quad\mbox{for}\quad\tau\in I.
\end{equation}
\end{theorem}
Theorem~\ref{Theorem:1.2} is closely related to \cite[Theorem~88]{HLP}, 
where equivalent means were discussed. 

Next we characterize $F$-concavities 
closed under positive scalar multiplication 
and positive exponentiation, respectively.
\begin{theorem}
\label{Theorem:1.3}
Let $\Omega$ be a convex domain in ${\bf R}^N$. 
Let $I$ be an open interval and let $F:I\to {\bf R}$ satisfy the following condition: 
\begin{itemize}
  \item[{\rm (F)}] 
$\quad F\in C^1(I),\quad F'>0\quad\mbox{on}\quad I.$
 \end{itemize}
Assume that $F$-concavity is nontrivial in $\Omega$ and possesses the following property:
\begin{itemize}
  \item[{\rm (P)}]  
  For any $f\in{\mathcal C}_\Omega[F]\cap{\mathcal A}^*_\Omega(I)$, 
  there exists $\epsilon\in(0,1)$ such that 
  $$
  \lambda f\in{\mathcal C}_\Omega[F]
  \quad\mbox{for}\quad \lambda\in(1-\epsilon,1+\epsilon).
  $$ 
\end{itemize}
Then the following assertions hold.
\begin{itemize}
  \item[{\rm (a)}] 
  If $I\subset(0,\infty)$, then there exists a triple $(\alpha,A_1,B_1)\in {\bf R}\times (0,\infty)\times {\bf R}$ such that
  $$
  F(\tau)=A_1\Phi_\alpha(\tau)+B_1\quad\mbox{for}\quad\tau\in I.
  $$
  \item[{\rm (b)}]
  If $I\subset(-\infty,0)$, then there exists a triple $(\beta,A_2,B_2)\in {\bf R}\times (0,\infty)\times {\bf R}$ such that
  $$
  F(\tau)=-A_2\Phi_\beta(-\tau)+B_2\quad\mbox{for}\quad\tau\in I.
  $$
  \item[{\rm (c)}] 
  If $0\in I$, then there exists a pair $(A_3,B_3)\in (0,\infty)\times {\bf R}$ such that 
  $$
  F(\tau)=A_3\tau+B_3\quad\mbox{for}\quad\tau\in I. 
  $$
\end{itemize}
\end{theorem}
\begin{theorem}
\label{Theorem:1.4}
Let $\Omega$ be a convex domain in ${\bf R}^N$. 
Let $I$ be an open interval in~$(0,\infty)$ and let $F:I\to {\bf R}$ satisfy condition~{\rm (F)}. 
Assume that $F$-concavity is nontrivial in $\Omega$ and possesses the following property: 
\begin{itemize}
  \item[{\rm (P')}]  
  For any $f\in{\mathcal C}_\Omega[F]\cap{\mathcal A}^*_\Omega(I)$, 
  there exists $\epsilon \in (0,1)$ such that 
  $$
  f^r\in{\mathcal C}_\Omega[F]
  \quad\mbox{for}\quad r\in(1-\epsilon,1+\epsilon).
  $$ 
\end{itemize}
Then the following assertions hold.
\begin{itemize}
  \item[{\rm (a')}] 
  If $I\subset(1,\infty)$, then there exists a triple $(\alpha,A_1,B_1)\in {\bf R}\times (0,\infty)\times {\bf R}$ such that
  $$
  F(\tau)=A_1\Phi_\alpha(\log\tau)+B_1\quad\mbox{for}\quad\tau\in I.
  $$
  \item[{\rm (b')}]
  If $I\subset(0,1)$, then there exists a triple $(\beta,A_2,B_2)\in {\bf R}\times (0,\infty)\times {\bf R}$ such that
  $$
  F(\tau)=-A_2\Phi_\beta(-\log\tau)+B_2\quad\mbox{for}\quad\tau\in I.
  $$
  \item[{\rm (c')}] 
  If $1\in I$, then there exists a pair $(A_3,B_3)\in (0,\infty)\times {\bf R}$ such that 
  $$
  F(\tau)=A_3\log\tau+B_3\quad\mbox{for}\quad\tau\in I. 
  $$
\end{itemize}
\end{theorem}
Theorem~\ref{Theorem:1.4} is a byproduct of Theorem~\ref{Theorem:1.3}. 
Theorems~\ref{Theorem:1.3} and \ref{Theorem:1.4} imply that, 
among nontrivial $F$-concavities satisfying condition~(F), the following hold:
\begin{itemize}
  \item[(i)] 
  only power concavity is closed under positive scalar multiplication;
  \item[(ii)] 
  only power log-concavity is closed under positive exponentiation;
  \item[(iii)] 
  {\em only log-concavity is closed under both positive scalar multiplication and positive exponentiation}.
\end{itemize}
In Section~3 
we also identify $F$-concavities closed under translation (see Theorem \ref{Theorem:3.1}).
\bigskip

Next we discuss the preservation of $F$-concavity by the Dirichlet heat flow. 
For any $\varphi\in L^\infty(\Omega)$, 
we denote by $e^{t\Delta_\Omega}\varphi$ the unique bounded solution to the problem
\begin{equation}
\label{eq:1.2}
\left\{
\begin{array}{ll}
\partial_t u=\Delta u & \mbox{in}\quad\Omega\times(0,\infty),\vspace{3pt}\\
u=0 & \mbox{on}\quad\partial\Omega\times(0,\infty)\quad\mbox{if}\quad\partial\Omega\not=\emptyset,\vspace{3pt}\\
u(\cdot,0)=\varphi(\cdot) & \mbox{in}\quad\Omega.
\end{array}
\right.
\end{equation}
We say that {\it $F$-concavity is preserved by the Dirichlet heat flow in $\Omega$} if 
$$
e^{t\Delta_\Omega}\varphi\in{\mathcal C}_\Omega[F]\quad\mbox{for}\quad t>0\quad\mbox{if}\quad 
\varphi\in{\mathcal C}_\Omega[F]\cap L^\infty(\Omega).
$$
In what follows, we will consider only {\em nonnegative initial data}, i.e. $\varphi\geq0$ in $\Omega$ (whence {\em nonnegative solutions}).

As mentioned in assertion~(C'),  log-concavity is preserved by  the Dirichlet heat flow.
In this connection one of the main motivations of this paper is the following natural question:
\begin{itemize}
  \item[(Q)] 
  What is the strongest (nontrivial) $F$-concavity preserved by the Dirichlet heat flow?
\end{itemize}
The answer may be not unique 
since it strongly depends on the interval $I$ where $F$ is defined. See Definition \ref{Definition:1.1}. 
When we consider nonnegative initial data, we can restrict our attention 
to the intervals~$I$ such that $\mbox{Int}\,I=(0,a)$ for some $a\in(0,\infty]$. 
When $I=[0,\infty)$, the answer we give here is {\em ``\,log-concavity"}: 
this is the content of the following Theorem~\ref{Theorem:1.5} (jointly with assertion~(C')). 
\begin{theorem}
\label{Theorem:1.5}
Let $\Omega$ be a convex domain in ${\bf R}^N$.
Let $F$ be admissible on $[0,\infty)$ such that $F(0)=-\infty$. 
If $F$-concavity is preserved by the Dirichlet heat flow in $\Omega$, 
then $F$-concavity is weaker than log-concavity in $\Omega$, that is, 
${\mathcal C}_\Omega[\Phi_0]\subset{\mathcal C}_\Omega[F]$. 
\end{theorem}

\begin{remark}
\label{Remark:1.3}
Let $F$ be admissible on $[0,\infty)$ and $\Omega={\bf R}^N$. 
If $F(0)>-\infty$, then Theorem~{\rm\ref{Theorem:1.1}~(b)} implies that $F$-concavity is trivial, 
that is, ${\mathcal C}_\Omega[F]$ consists only of nonnegative constant functions in ${\bf R}^N$ 
and $F$-concavity is trivially preserved by the Dirichlet heat flow. 
\end{remark}

We already dealt with question (Q) in our previous paper \cite{IST}*{Theorem 3.2} and our answer was different: {\em $(1/2)$-log-concavity}\,!
Beyond appearances, Theorem \ref{Theorem:1.5} and the conclusion of \cite{IST} are not in contrast. 
Indeed, the intervals considered are different 
(and even the definition of $F$-concavity is different in \cite{IST}). 
For a further discussion, see  Section \ref{Section:1.5}. 

\vspace{3pt}

The rest of this paper is organized as follows. 
In Section~2 we prove Theorems~\ref{Theorem:1.1} and \ref{Theorem:1.2}. 
In Section~3 we prove Theorem~\ref{Theorem:1.3}.  
As applications of Theorem~\ref{Theorem:1.3}, 
we prove Theorem~\ref{Theorem:1.4} and give a characterization 
of $F$-concavities closed under translation (see Theorem \ref{Theorem:3.1}).
In Section~4 we prove Theorem~\ref{Theorem:1.5} and discuss in detail 
the relation between log-concavity and $(1/2)$-log-concavity.
In Section~5 we present some open problems related to the Dirichlet heat flow.
%%%%%%%%%%%%%%%%%%%%%%%%%%%%
\section{Proofs of Theorems~\ref{Theorem:1.1} and \ref{Theorem:1.2}}
%%%%%%%%%%%%%%%%%%%%%%%%%%%%
We prove Theorems~\ref{Theorem:1.1} and \ref{Theorem:1.2} to characterize nontrivial $F$-concavities. 
\vspace{5pt}
\newline
{\bf Proof of Theorem~\ref{Theorem:1.1}.}
Let $\Omega$ be a convex domain with $\Omega\not={\bf R}^N$. 
Then we find $x_*\in\partial\Omega$ and $\nu\in{\bf R}^N$ such that 
$$
\Omega\subset\{x\in{\bf R}^N\,:\,\langle x-x_*,\nu\rangle> 0\}
$$
(see e.g. \cite[Theorem 1.3.2]{Sh}). 
Here $\langle \cdot,\cdot\rangle$ is the standard inner product in ${\bf R}^N$. 
Let $F$ be admissible on $I$. 
Set 
$$
g(x):=F(b)-(F(b)-F(a))e^{-\langle x-x_*,\nu\rangle}\quad\mbox{for}\quad x\in\Omega,
$$
where $-\infty<a<b<\infty$ and $[a,b]\subset\mbox{Int}\,I$. 
Then $g$ is strictly concave in $\Omega$ and $F(a)<g(x)<F(b)$ in $\Omega$.
Since $F$ is strictly increasing and continuous on $[a,b]$, 
the inverse function $G$ of $F$ can be defined on $[F(a),F(b)]$ 
and set $f(x):=G(g(x))\in[a,b]$ for $x\in\Omega$. 
Then we easily see that $f$ is a nonconstant $F$-concave function in $\Omega$, 
that is, $F$-concavity is nontrivial in $\Omega$. 
Thus assertion~(a) follows. 

We prove assertion~(b). 
Let $\Omega={\bf R}^N$ and $\inf_I F=-\infty$. 
Let $a$, $b$ be as in the above. 
Assume that $\inf I\not\in I$. 
Since $F$ is admissible on $I$ and $\inf_I F=-\infty$, 
we see that 
$$
\lim_{\tau\to \inf I+0}F(\tau)=-\infty,
$$ 
which implies that the inverse function $G$ of $F$ can be defined on $(-\infty,F(b)]$. 
Set 
$$
\tilde{g}(x):=F(b)-e^{-x_1},\qquad 
\tilde{f}(x):=G(\tilde{g}(x)),
$$
for $x=(x_1,x')\in{\bf R}\times{\bf R}^{N-1}$.
Then $\tilde{f}$ is a nonconstant $F$-concave function in $\Omega$. 
If $\inf I\in I$, then $F(\inf I)=-\infty$ and set 
$$
\hat{g}(x):=
\left\{
\begin{array}{ll}
F(b)-(F(b)-F(a))e^{-x_1} & \mbox{for}\quad x_1>0,\vspace{3pt}\\
-\infty & \mbox{for}\quad x_1\le 0,
\end{array}
\right.
\qquad 
\hat{f}(x):=
\left\{
\begin{array}{ll}
G(\hat{g}(x)) & \mbox{for}\quad x_1>0,\vspace{3pt}\\
\inf I & \mbox{for}\quad x_1\le 0,
\end{array}
\right.
$$
for $x=(x_1,x')\in{\bf R}\times{\bf R}^{N-1}$. 
Then $\hat{f}$ is a nonconstant $F$-concave function in $\Omega$. 
Therefore we see that $F$-concavity is nontrivial in ${\bf R}^N$ if $\inf_I F=-\infty$. 

On the other hand, let $f\in{\mathcal C}_{{\bf R}^N}[F]$ be such that $F(f)$ is  bounded from below. Since a nonconstant concave function in ${\bf R}^N$ must be unbounded from below, 
we see that $F(f)$ is constant, 
whence $f$ must be constant and
 assertion~(b) follows. 
$\Box$\vspace{5pt}
\newline
{\bf Proof of Theorem~\ref{Theorem:1.2}.}
If \eqref{eq:1.4} holds, then it easily follows that ${\mathcal C}_\Omega[F_1]={\mathcal C}_\Omega[F_2]$. 
So we have only to prove that ${\mathcal C}_\Omega[F_1]={\mathcal C}_\Omega[F_2]$ implies \eqref{eq:1.4}. 
The proof is by contradiction. 
Assume that \eqref{eq:1.4} does not hold. 
Then we can assume, without loss of generality,  that 
\begin{equation}
\label{eq:2.1}
F_1(a)=F_2(a)=0,\qquad F_1(b)=F_2(b)=1,\qquad
F_1(c)\not=F_2(c),
\end{equation}
holds for some $a$, $b$, $c\in I$ with $a<c<b$ (see also Remark \ref{rem-F}). 
The strict monotonicity of $F_1$ and $F_2$ implies that $F_1(c)$, $F_2(c)\in(0,1)$. 

Let $\Omega$ be a convex domain in ${\bf R}^N$ with $\Omega\not={\bf R}^N$. 
Let $a'\in I$ be such that $a'<a$. 
Thanks to \eqref{eq:1.1}, 
we can assume, without loss of generality, that  $0\in \Omega$ and
$$
\{ s e_1\,:\,F_1(a')<s\le 1\}\subset\Omega
\subset\{(x_1,x')\,:\,x_1>F_1(a'),\,\,x'\in{\bf R}^{N-1}\},
$$ 
where  $e_1=(1,0,\ldots,0 )\in {\bf R}^N$.
Set 
$$
g(x)=\min\{x_1,1\}\quad\mbox{for}\quad x=(x_1,x')\in\Omega.
$$ 
Then $g$ is a concave function in $\Omega$ such that 
$\{g(x)\,:\,x\in\Omega\}=(F_1(a'),1]$. 
Since the inverse function $G_1$ of $F_1$ can be defined on $(F_1(a'),1]$, 
we set 
$$
f(x):=G_1(g(x))\quad\mbox{for}\quad x\in\Omega.
$$ 
Then $f$ is $F_1$-concave in $\Omega$. 
Since ${\mathcal C}_\Omega[F_1]={\mathcal C}_\Omega[F_2]$, 
we see that 
the function $h=F_2(f)$ is concave in $\Omega$. Furthermore, 
$$
h(0)=F_2(G_1(0))=0,\qquad h(e_1)=F_2(G_1(1))=1.
$$
These together with the concavity of $w$ imply that 
$$
h(s e_1)\ge s,\quad\mbox{that is},\quad
f(s e_1)\ge G_2(s),
$$
for $s\in (0,1]$, where $G_2$ is the inverse function of $F_2$ on $[0,1]$.
Then we deduce that 
$$
c=G_1(F_1(c))=G_1(g(F_1(c)e_1))=f(F_1(c))\ge G_2(F_1(c)),
$$
that is, $F_2(c)\ge F_1(c)$. 
Similarly, we can show that $F_2(c)\le F_1(c)$. Then $F_2(c)=F_1(c)$, 
which contradicts \eqref{eq:2.1}. 
Thus Theorem~\ref{Theorem:1.2} follows in the case $\Omega\not={\bf R}^N$.

Let $\Omega={\bf R}^N$. Since $F_1$-concavity is nontrivial, 
by Theorem~\ref{Theorem:1.1}~(b) we see that $F_1$ and $F_2$ are unbounded below. 
Then we apply the same argument as in the case $\Omega\not={\bf R}^N$ 
to find a contradiction. Thus Theorem~\ref{Theorem:1.2} follows in the case $\Omega={\bf R}^N$. 
The proof is complete.
$\Box$
%%%%%%%%%%%%%%%%%%%%%%%%%%%%
\section{Proofs of Theorems~\ref{Theorem:1.3} and \ref{Theorem:1.4}}
%%%%%%%%%%%%%%%%%%%%%%%%%%%%
In this section we characterize $F$-concavities closed 
under positive scalar multiplication, under positive exponentiation and under translation.
By the letter $C$
we denote a generic positive constant and this may have different values even within the same line. 
\vspace{5pt}
\newline
{\bf Proof of Theorem~\ref{Theorem:1.3}.}
Consider the case $\Omega\not={\bf R}^N$. 
Let $x_*\in\partial\Omega$ and $\nu\in{\bf R}^N$ be as in the proof of Theorem \ref{Theorem:1.1}. 
We can assume, without loss of generality, that $x_*=0$ and $\nu=e_1$,
whence  $\Omega\subset {\bf R}^N_+:=\{(x_1,x')\in{\bf R}^N\,:\,x_1>0,\,\,x'\in{\bf R}^{N-1}\}$ 
and $0\in\partial\Omega$.

Let $F$ be admissible on an open interval $I$. 
Let $-\infty<m_1<m_2<m_2'<\infty$ be such that $[m_1,m_2']\subset I$. 
Set 
$$
d:=\min\{1\,;\sup\{x_1\,:\, (x_1,0)\in\Omega\}\}>0,
\quad
S:=(0,d),\quad
k:=\frac{F(m_2)-F(m_1)}{d}>0.
$$
Let $g=g(s)$ be a smooth concave function on $(0,\infty)$ such that 
\begin{equation}
\label{eq:3.1}
g(s)=F(m_1)+ks\quad\mbox{for}\quad s\in S,
\quad
F(m_1)<g(s)<F(m_2')\quad\mbox{for}\quad s\in(0,\infty).
\end{equation}
Then 
$\{g(s)\,:\,s\in S\}=(F(m_1),F(m_2))$.
Since the inverse function $G$ of $F$ can be defined on the interval $(F(m_1),F(m_2'))$,  
we set 
$$f(s):=G(g(s))\quad\text{for }s>0\,.
$$
It follows from $G\in C^1((F(m_1),F(m_2)])$ that $f\in C^1(S)$. 
Furthermore, we have
\begin{equation}
\label{eq:3.2}
\begin{split}
 & \frac{d}{ds}F(f(s))=F'(f(s))\frac{d}{ds}f(s)=k\quad\mbox{and}\quad
\frac{d}{ds}f(s)>0\quad
\mbox{for}\quad s\in S,\\
 & \{f(s)\,:\,s\in S\}=(m_1,m_2)\quad\mbox{and}\quad
  \{f(s)\,:\,s>0\}\subset (m_1,m_2').
\end{split}
\end{equation}
Set 
$$
\tilde{f}(x):=f(x_1)\quad\mbox{for}\quad x=(x_1,x')\in\Omega. 
$$
Then $\tilde{f}$ is $F$-concave in $\Omega$. 
Since $\tilde{f}\in{\mathcal C}_\Omega[F]\cap{\mathcal A}^*_\Omega(I)$, 
by property~(P) we find $\epsilon\in(0,1)$ such that 
$$
\lambda\tilde{f}\in{\mathcal C}_\Omega[F]\quad\mbox{for}\quad\lambda\in(1-\epsilon,1+\epsilon).
$$
This implies that 
\begin{equation*}
\begin{split}
{\mathcal F}(x,y,\lambda,\mu): & =
F(\lambda \tilde{f}((1-\mu)x+\mu y))-(1-\mu)F(\lambda \tilde{f}(x))-\mu F(\lambda \tilde{f}(y))\\
 & =\,F(\lambda f((1-\mu)x_1+\mu y_1))-(1-\mu)F(\lambda f(x_1))-\mu F(\lambda f(y_1))\ge 0
\end{split}
\end{equation*}
for $x=(x_1,x')$, $y=(y_1,y')\in\Omega$, $\mu\in[0,1]$ and $\lambda\in(1-\epsilon,1+\epsilon)$. 
On the other hand, by \eqref{eq:3.1} we see that 
$$
{\mathcal F}(x,y,1,\mu)=g((1-\mu)x_1+\mu y_1)-(1-\mu)g(x_1)-\mu g(y_1)=0
$$
for $x=(x_1,x')$, $y=(y_1,y')\in \Omega$ with $x_1$, $y_1\in S$ and $\mu\in[0,1]$. 
These imply that 
\begin{equation*}
\begin{split}
0 & =\frac{\partial}{\partial\lambda}{\mathcal F}(x,y,\lambda,\mu)\biggr|_{\lambda=1}\\
 & =F'(f((1-\mu)x_1+\mu y_1))f((1-\mu)x_1+\mu y_1)-(1-\mu)F'(f(x_1))f(x_1)-\mu F'(f(y_1))f(y_1)
\end{split}
\end{equation*}
for $x=(x_1,x')$, $y=(y_1,y')\in \Omega$ with $x_1$, $y_1\in S$ and $\mu\in[0,1]$, 
that is, 
$$
F'(f((1-\mu)s_1+\mu s_2))f((1-\mu)s_1+\mu s_2)=(1-\mu)F'(f(s_1))f(s_1)+\mu F'(f(s_2))f(s_2)
$$
for $s_1$, $s_2\in S$. 
Then we see that the function $F'(f(s))f(s)$ is a linear function on $S$
and find $a$, $b\in{\bf R}$ such that 
\begin{equation}
\label{eq:3.3}
F'(f(s))f(s)=as+b\quad\mbox{for}\quad s\in S. 
\end{equation} 
This together with \eqref{eq:3.2} implies that 
\begin{equation}
\label{eq:3.4}
f(s)\biggr/\frac{d}{ds}f(s)=\frac{as+b}{k}\quad\mbox{for}\quad s\in S. 
\end{equation} 

We consider the case $0\notin I$ and set 
\[
\sigma_I:=\begin{cases} 1 & \text{if\ } I\subset(0,\infty),\\ -1 & \text{if\ } I\subset(-\infty,0).\end{cases}
\]
We  prove assertions~(a) and (b).
Since $ f(s)\in I$ for $s\in S$, 
by \eqref{eq:3.3} we see that
$$
  F'(f(s))\cdot \sigma_I f(s)=\sigma_I (as+b)>0
$$ 
for $s\in S$, 
in particular, $\sigma_Ib>0$ if $a=0$.
It follows from \eqref{eq:3.4} that
$$
\frac{d}{ds}\log |f(s)|=\frac{1}{f(s)}\frac{d}{ds}f(s)=\frac{k}{as+b},
$$
that is, 
$$
\log|f(s)|=
\left\{
\begin{array}{ll}
\displaystyle{\frac{k}{a}}\log|as+b|+C & \mbox{if}\quad a\not=0,\vspace{7pt}\\
\displaystyle{\frac{k}{b}s}+C & \mbox{if}\quad a=0,
\end{array}
\right.
$$
for $s\in S$. Then 
$$
f(s)=
\left\{
\begin{array}{ll}
\sigma_I e^C(|as+b|)^{\frac{k}{a}} & \mbox{if}\quad a\not=0,\vspace{7pt}\\
\sigma_I e^C\displaystyle{\exp\left(\frac{k}{b}s\right)} & \mbox{if}\quad a=0,
\end{array}
\right.
\qquad
\frac{d}{ds}f(s)= \frac{k}{as+b} f(s),
$$
for $s\in S$. 
By \eqref{eq:3.2} we obtain
\begin{equation}
\label{eq:3.5}
\begin{split}
F'(f(s))
 & =k\biggr/\frac{d}{ds}f(s)
   =\frac{as+b}{f(s)}
 =
\left\{
\begin{array}{ll}
e^{-\frac{aC}{k}} |f(s)|^{-1+\frac{a}{k}} & \mbox{if}\quad a\not=0,\vspace{7pt}\\
b f(s)^{-1} & \mbox{if}\quad a=0,
\end{array}
\right.
\end{split}
\end{equation}
for $s\in S$, which implies that
$$
F'(\tau)
=
\left\{
\begin{array}{ll}
e^{-\frac{aC}{k}} |\tau|^{-1+\frac{a}{k}} & \mbox{if}\quad a\not=0,\vspace{7pt}\\
b \tau^{-1} & \mbox{if}\quad a=0,
\end{array}
\right.
$$
for $\tau\in(m_1,m_2)$. 
Integrating this equality yields 
$$
F(\tau)=A \Phi_\beta( \sigma_I \tau)+C
$$
for $\tau\in(m_1,m_2)$, 
where
\[
\beta=\frac{a}{k}\in{\bf R}, \quad
A=\begin{cases}
\displaystyle  \sigma_I  e^{-\frac{aC}{k}} & \mbox{if}\quad a\not=0, \vspace{7pt}\\
b  & \mbox{if}\quad a=0.
\end{cases}
\]
Notice that  $\sigma_I  A>0$.
Since $m_1$ and $m_2$ are arbitrary, we obtain $F( \sigma_I \tau)=A  \Phi_\beta(\tau)+C$
for $\tau\in I$.
Thus assertions~(a) and (b) follow.

We consider the case $0\in I$ and prove assertion~(c). 
We can assume, without loss of generality, that $m_1<0<m_2$. 
Since there exists $s_*\in S$ such that $f(s_*)=0$,  
we find $a\in{\bf R}$ such that
$$
F'(f(s))f(s)=a(s-s_*)\quad\mbox{for}\quad s\in S,
$$
instead of \eqref{eq:3.3}. 
Since $F'>0$ and $f$ is strictly increasing on $S$, 
we see that $a>0$.
Furthermore, 
\begin{equation}
\label{eq:3.8}
f(s)<0\quad\mbox{for}\quad s\in(0,s_*)\quad\mbox{and}\quad
f(s)>0\quad\mbox{for}\quad t\in(s_*,d).
\end{equation}
Repeating the above argument, we see that 
$$
\frac{d}{ds}\log |f(s)|=\frac{1}{f(s)}\frac{d}{ds}f(s)=\frac{k}{a(s-s_*)}
\quad\mbox{for}\quad s\in S\setminus\{s_*\}.
$$
Then we find $C_1$, $C_2\in{\bf R}$ such that 
$$
\log|f(s)|=
\left\{
\begin{array}{ll}
\displaystyle{\frac{k}{a}}\log|a(s-s_*)|+C_1 & \quad\mbox{for}\quad 0<s<s_*,\vspace{7pt}\\
\displaystyle{\frac{k}{a}}\log|a(s-s_*)|+C_2 & \quad\mbox{for}\quad s_*<s<d,
\end{array}
\right.
$$
which together with \eqref{eq:3.8} implies that 
$$
f(s)=-e^{C_1}a(s_*-s)^{\frac{k}{a}}\quad\mbox{for}\quad s\in(0,s_*),
\quad
f(s)=e^{C_2}a(s-s_*)^{\frac{k}{a}}\quad\mbox{for}\quad s\in(s_*,d).
$$
Recalling $f\in C^1(S)$ and $a>0$, 
by \eqref{eq:3.2} 
we see that $k/a=1$ and $C_1=C_2$, that is, 
$$
f(s)=e^{C_1}(s-s_*)
$$
for $s\in S=(0,d)$. 
Similarly to \eqref{eq:3.5}, 
we have 
$$
F'(f(s))=k\biggr/\frac{d}{ds}f(s)=ke^{-C_1}
$$
for $s\in S$, which implies that 
$F'(\tau)=ke^{-C_1}$
for $\tau \in (m_1, m_2)$.
Then, repeating the above argument, 
we deduce that 
$$
F(\tau)=A \tau+C
$$
for $\tau\in I$, where $A=ke^{-C_1}>0$. 
Thus assertion~(c) follows. 
Therefore the proof of Theorem~\ref{Theorem:1.3} is complete in the case $\Omega\not={\bf R}^N$. 

It remains to consider the case $\Omega={\bf R}^N$. 
Since $F$-concavity is nontrivial, by Theorem~\ref{Theorem:1.1}~(b) we see that $\inf_I F=-\infty$. 
Let $-\infty<m_1<m_2<m_2'<\infty$ be such that $[m_1,m_2']\subset I$. 
Let $\tilde{g}=\tilde{g}(s)$ be a smooth concave function on ${\bf R}$ such that 
$$
\tilde{g}(s)=F(m_1)+ks\quad\mbox{for}\quad s\in \tilde{S}:=(-\infty,1),
\quad
F(m_1)<\tilde{g}(s)<F(m_2')\quad\mbox{for}\quad s\in(-\infty,\infty).
$$
Since the inverse function $G$ of $F$ can be defined on $(-\infty,F(m_2'))$, 
we set $\tilde{f}(s)=G(\tilde{g}(s))$ for $s\in\tilde{S}$. 
Then, repeating the same argument as in the case $\Omega\not={\bf R}^N$, 
we obtain assertions~(a), (b) and (c). 
Thus Theorem~\ref{Theorem:1.3} follows.
$\Box$\vspace{5pt}

Next we apply Theorem~\ref{Theorem:1.3} to prove Theorem~\ref{Theorem:1.4}. 
\vspace{3pt}
\newline
{\bf Proof of Theorem~\ref{Theorem:1.4}.}
Let $F$ be admissible on an open interval $I\subset(0,\infty)$. 
Set 
$$
J:=\{\log\tau\,:\,\tau\in I\},
\qquad
\tilde{F}(t):=F\left(e^t\right)\quad\mbox{for}\quad t\in J. 
$$
Then $J$ is an open interval and $\tilde{F}$ is admissible on $J$. 
Furthermore, 
$$
f\in{\mathcal C}_\Omega[F]
\quad\mbox{if and only if}\quad 
\log f\in{\mathcal C}_\Omega[\tilde{F}].
$$
It follows from property~(P') that 
$\tilde{F}$ possesses property~(P). 
Therefore we find that 
assertions~(a), (b) and (c) in Theorem~\ref{Theorem:1.3} 
hold with $F$ and $I$ replaced by $\tilde{F}$ and $J$, respectively. 
This implies the desired conclusion, and Theorem~\ref{Theorem:1.4} follows.
$\Box$
\vspace{5pt}

Similarly to the proof of Theorem~\ref{Theorem:1.4}, 
we give a characterization of $F$-concavities closed under translation.
\begin{theorem}
\label{Theorem:3.1}
Let $\Omega$ be a convex domain in ${\bf R}^N$. 
Let $I$ be an open interval and let $F:I\to {\bf R}$ satisfy condition~{\rm (F)}. 
Assume that $F$-concavity is nontrivial in $\Omega$ and it possesses property:
\begin{itemize}
  \item[{\rm (P'')}]  
  For any $f\in{\mathcal C}_\Omega[F]\cap{\mathcal A}^*_\Omega(I)$, 
  there exists $\epsilon>0$ such that 
  $$
  f+c\in{\mathcal C}_\Omega[F] \quad\mbox{for}\quad  c\in(-\epsilon,\epsilon).
  $$ 
\end{itemize}
Then there exists a triple $(\alpha,A,B)\in {\bf R}\times (0,\infty)\times {\bf R}$ such that
$$
F(\tau)=A\Phi_{\alpha}(e^\tau)+B\quad\mbox{for}\quad \tau\in I.
$$
\end{theorem}
{\bf Proof.}
Let $F$ be admissible on an open interval $I$. 
Set 
$$
J:=\{e^\tau\,:\,\tau\in I\}\subset(0,\infty),
\qquad
\tilde{F}(t):=F(\log t)\quad\mbox{for}\quad t\in J. 
$$
Then $J$ is an open interval and $\tilde{F}$ is admissible on $J$. 
Furthermore, 
$$
f\in{\mathcal C}_\Omega[F]
\quad\mbox{if and only if}\quad 
e^f\in{\mathcal C}_\Omega[\tilde{F}].
$$
It follows from property~(P'') that 
$\tilde{F}$ possesses property~(P).
Therefore, by Theorem~\ref{Theorem:1.3}~(a)
we find a triple $(\alpha,A,B)\in{\bf R}\times(0,\infty)\times{\bf R}$ such that 
$$
\tilde{F}(\tau)=A\Phi_{\alpha}(\tau)+B
$$
for $\tau \in  J$. This implies the desired conclusion.
$\Box$
%%%%%%%%%%%%%%%%%%%%%%%%%%%%%%%%%%%%
%%%%%%%%%%%%%%%%%%%%%%%%%%%%%%%%%%%%
\section{$F$-concavity and Dirichlet heat flow}
%%%%%%%%%%%%%%%%%%%%%%%%%%%%%%%%%%%%
%%%%%%%%%%%%%%%%%%%%%%%%%%%%%%%%%%%%
We prove Theorem~\ref{Theorem:1.5} and discuss the relation between log-concavity and $(1/2)$-log-concavity. 
%%%%%%%%%%%%%%%%%%%
\subsection{Proof of Theorem~\ref{Theorem:1.5}}
%%%%%%%%%%%%%%%%%%%
The following lemma is crucial for the proof of Theorem~\ref{Theorem:1.5}. 
\begin{lemma}
\label{Lemma:4.1}
Let $\Omega$ be a convex domain in ${\bf R}^N$. 
Let $F$ be admissible on $[0,\infty)$ such that $F(0)=-\infty$.
Assume that $F$-concavity is preserved by the Dirichlet heat flow in $\Omega$. 
Then the function 
$$
[0,\infty)\ni s\mapsto F\left(ke^{-s^2}\right)
$$ 
is concave for every $k>0$. 
\end{lemma}
{\bf Proof.}
We modify the proof of \cite[Theorem~3.2]{IST} to prove Lemma~\ref{Lemma:4.1}. 
The proof is by contradiction. 
We assume that the function 
$$
[0,\infty)\ni s\mapsto F\left(ke^{-s^2}\right)
$$ 
is not concave for some $k>0$. 
Then we find $s_1$, $s_2>0$ and $\mu\in(0,1)$ such that 
\begin{equation}
\label{eq:4.1}
F\left(ke^{-\{(1-\mu)s_1+\mu s_2\}^2}\right)<(1-\mu)F\left(ke^{-s_1^2}\right)+\mu F\left(ke^{-s_2^2}\right). 
\end{equation}

Let 
$$
u(x,t):=\left[e^{t\Delta_{{\bf R}^N}}\chi_{B(0,1)}\right](x)
=(4\pi t)^{-\frac{N}{2}}\int_{B(0,1)}e^{-\frac{|x-y|^2}{4t}}\,dy
$$ 
for $x\in{\bf R}^N$ and $t>0$. 
Then, for any $L>0$,  
\begin{equation}
\label{eq:4.2}
\lim_{t\to\infty}\sup_{\{|x|\le L\sqrt{t}\}}\left|(4\pi t)^{\frac{N}{2}}|B(0,1)|^{-1}u(x,t)-e^{-\frac{|x|^2}{4t}}\right|=0,
\end{equation}
where $|B(0,1)|$ is the volume of $B(0,1)$. 
Set 
$$
x(t):=2\sqrt{t} s_1 e_1,
\qquad
y(t):=2\sqrt{t} s_2 e_1,
$$
for $t>0$.
Taking a sufficiently large $T>0$, 
we observe from \eqref{eq:4.1} and \eqref{eq:4.2} that 
\begin{equation*}
\begin{split}
 & F\left(k(4\pi T)^{\frac{N}{2}}|B(0,1)|^{-1}u((1-\mu)x(T)+\mu y(T),T)\right)\\
 & <(1-\mu)F\left(k(4\pi T)^{\frac{N}{2}}|B(0,1)|^{-1}u(x(T),T)\right)
+\mu F\left(k(4\pi T)^{\frac{N}{2}}|B(0,1)|^{-1}u(y(T),T)\right).
\end{split}
\end{equation*}
Setting
$$
v(x,t):=k(4\pi T)^{\frac{N}{2}}|B(0,1)|^{-1}u(x,t)
=k(4\pi T)^{\frac{N}{2}}|B(0,1)|^{-1}\left[e^{t\Delta_{{\bf R}^N}}\chi_{B(0,1)}\right](x),
$$
we obtain 
\begin{equation}
\label{eq:4.3}
F\left(v((1-\mu)x(T)+\mu y(T),T)\right)<(1-\mu)F\left(v(x(T),T)\right)+\mu F\left(v(y(T),T)\right).
\end{equation}
This implies that $v(\cdot,T)$ is not $F$-concave in ${\bf R}^N$. 

Let $\Omega$ be a convex domain in ${\bf R}^N$. 
We can assume, without loss of generality, that $0\in\Omega$. 
Set $\Omega_n=n\Omega$ for $n=1,2,\dots$. 
For sufficiently large $n$, let 
$$
v_n(x,t):=k(4\pi T)^{\frac{N}{2}}|B(0,1)|^{-1}\left[e^{t\Delta_{\Omega_n}}\chi_{B(0,1)}\right](x)
$$
for $x\in\Omega_n$ and $t>0$. 
Since $\Omega_n\to{\bf R}^N$ as $n\to\infty$, 
we see that 
$$
\lim_{n\to\infty}v_n(x,t)=v(x,t)\quad\mbox{for}\quad (x,t)\in{\bf R}^N\times(0,\infty).
$$
This together with \eqref{eq:4.3} yields
\begin{equation}
\label{eq:4.4}
F\left(v_n((1-\mu)x(T)+\mu y(T),T)\right)<(1-\mu)F\left(v_n(x(T),T)\right)+\mu F\left(v_n(y(T),T)\right)
\end{equation}
for sufficiently large $n$. 

For sufficiently large $n$, set $w_n(x,t):=v_n(nx,n^2t)$ for $x\in\Omega$ and $t>0$. 
Then $w_n$ satisfies 
$$
\left\{
\begin{array}{ll}
\partial_t w_n=\Delta w_n\quad & \mbox{in}\quad\Omega\times(0,\infty),\vspace{3pt}\\
w_n=0\quad & \mbox{in}\quad\partial\Omega\times(0,\infty),\vspace{3pt}\\
w_n(\cdot,0)=k(4\pi T)^{\frac{N}{2}}|B(0,1)|^{-1}\chi_{B(0,n^{-1})}\quad & \mbox{in}\quad\Omega.
\end{array}
\right.
$$
In addition,  $w_n(\cdot,0)$ is $F$-concave in $\Omega$ due to $F(0)=-\infty$. 
Furthermore, for sufficiently large~$n$, 
it follows from~\eqref{eq:4.4} that
\begin{equation*}
\begin{split}
 & F\left(w_n((1-\mu)n^{-1}x(T)+\mu n^{-1}y(T),n^{-2}T)\right)\\
 & <(1-\mu)F\left(w_n(n^{-1}x(T),n^{-2}T)\right)+\mu F\left(w_n(n^{-1}y(T),n^{-2}T)\right), 
\end{split}
\end{equation*}
which implies that $w_n(\cdot,n^{-2}T)$ is not $F$-concave in $\Omega$. 
These imply that $F$-concavity is not preserved by the Dirichlet heat flow in $\Omega$. 
This is a contradiction. 
Thus Lemma~\ref{Lemma:4.1} follows.~$\Box$\vspace{5pt}

Next we characterize log-concavity. 
Let $\rho\in C^\infty({\bf R})$ be even such that
$$
\rho\ge 0\quad\mbox{in}\quad{\bf R},
\qquad
\rho=0\quad\mbox{outside}\quad[-1,1],
\qquad
\int_{{\bf R}}\rho(s)\,ds=1.
$$
For $n=1,2,\dots$ set $\rho_n(s):=n\rho(ns)$ for $s\in{\bf R}$. Then 
$$
\rho_n\ge 0\quad\mbox{in}\quad{\bf R},
\qquad
\rho_n=0\quad\mbox{outside}\quad[-1/n,1/n],
\qquad
\int_{{\bf R}}\rho_n(s)\,ds=1.
$$
\begin{lemma}
\label{Lemma:4.2}
Let $F$ be admissible on $[0,\infty)$.
Then  the function 
$$
[0,\infty)\ni s\mapsto F\left(ke^{-s^2}\right)
$$ 
is concave for every $k>0$
if and only if the function $H(t):=F\left(e^t\right)$  is concave for $t\in{\bf R}$. 
Furthermore, in both cases,  $F(0)=-\infty$.
\end{lemma}
{\bf Proof.} 
First assume that $H$ is concave.
Let  $k>0$ and  set $\sigma:=\log k$.
Then, for $s_0,s_1 \geq 0$ and $\mu\in[0,1]$, 
we have
\begin{align*}\label{Fk-conc1}
F\left(ke^{-\{(1-\mu)s_0+\mu s_1\}^2}\right)
&=F\left(e^{\sigma-\{(1-\mu)s_0+\mu s_1\}^2}\right)\\
&\geq F\left(e^{\sigma-\{(1-\mu)s_0^2+\mu s_1^2\}}\right) 
=F\left(e^{(1-\mu)(\sigma-s_0^2)+\mu(\sigma- s_1^2)}\right)\\
&\geq(1-\mu)F\left(e^{\sigma-s_0^2}\right)+\mu F\left(e^{\sigma- s_1^2}\right)  
=(1-\mu)F\left(ke^{-s_0^2}\right)+\mu F\left(ke^{- s_1^2}\right),
\end{align*}
where the first inequality follows from  the monotonicity of $F$ together with the convexity of  
$[0,\infty)\ni s\mapsto s^2$.
Thus  we obtain the concavity of $[0,\infty)\ni s\mapsto F(ke^{-s^2})$.

Conversely,  given $t_0,\,t_1\in {\bf R}$ and $\mu\in[0,1]$, we prove that 
\begin{equation}\label{F-conc}
F\left(e^{(1-\mu)t_0+\mu t_1}\right)\geq (1-\mu)F\left(e^{t_0}\right)+\mu F\left(e^{t_1}\right).
\end{equation}
Setting 
\begin{equation*}
T:=\max\{t_0, t_1\},\quad
s_0=\sqrt{T-t_0},\quad 
s_1=\sqrt{T-t_1},\quad
k:=e^{T},
\end{equation*}
we see that
\begin{equation}\label{F-conc3} 
\begin{split}
(1-\mu)F\left( e^{t_0}\right)+\mu F\left(e^{t_1}\right)
&=(1-\mu)F\left(e^{T-s_0^2}\right)+\mu F\left(e^{T-s_1^2}\right) \\
&=(1-\mu)F\left(ke^{-s_0^2}\right)+\mu F\left(ke^{-s_1^2}\right)\\
&\leq F\left(ke^{-\{(1-\mu)s_0+\mu s_1\}^2}\right)
=F\left(e^{T-\left\{(1-\mu)\sqrt{T-t_0}+\mu\sqrt{T-t_1}\,\right\}^2}\right),
\end{split}
\end{equation}
where we used  the concavity of $F(ke^{-s^2})$ in the inequality.
Since
$$
\lim_{T\to+\infty}
\left[T-\left\{(1-\mu)\sqrt{T-t_0}+\mu\sqrt{T-t_1}\,\right\}^2\right]=(1-\mu)t_0+\mu t_1
$$
and $F$ is continuous in $(0,+\infty)$, passing to the limit as $T$ goes to $+\infty$ in \eqref{F-conc3} gives \eqref{F-conc}.

Finally, to prove that $F(0)=-\infty$ in both cases, 
just notice that 
the concavity of $H$ assures the existence of  $m\in{\bf R}$ with 
$$
H(t)\leq H(0)+mt \qquad  \text{for}\quad t\in {\bf R}.
$$
Since $H$ is strictly increasing, it must be $m>0$. Then
$$
\lim_{\tau\to+0}F(\tau)=
\lim_{t\to-\infty}F\left(e^t\right)=
\lim_{t\to-\infty}H(t)=-\infty\,,
$$
and the proof is complete.
$\Box$ 
\vspace{5pt}

Now we complete the proof of Theorem~\ref{Theorem:1.5}.
\vspace{5pt}
\newline
{\bf Proof of Theorem~\ref{Theorem:1.5}.}
Let $\Omega$ be a convex domain in ${\bf R}^N$. 
Assume that $F$-concavity is preserved by the Dirichlet heat flow in $\Omega$. 
It follows from Lemmas~\ref{Lemma:4.1} and~\ref{Lemma:4.2} that the function $H(t)=F\left(e^t\right)$ is concave for $t\in{\bf R}$.
To conclude the proof, we have to prove that $\mathcal{C}_\Omega[\Phi_0]\subseteq\mathcal{C}_\Omega[F]$, that is,  $F(f)$ is concave in $\Omega$ for every log-concave function $f$.

Let $f$ be a  log-concave function in $\Omega$. 
There exists a concave function $h:\Omega \to \mathbb{R}\cup\{-\infty\}$ such that $f=e^h$, that is,   
$F(f)=H(h)$ in $\Omega$.
Let $x,y\in\Omega$ with $f(x)f(y)>0$ and $\mu\in[0,1]$. 
The concavity of $h$ together with the monotonicity and the concavity of $H$ implies that 
$$
\begin{array}{rl}
F\left(f((1-\mu)x+\mu y)\right)&=H\left(h((1-\mu)x+\mu y)\right)\\
\\
&\ge H\left((1-\mu)h(x)+\mu h(y)\right)\\
\\
&\geq (1-\mu)H(h(x))+\mu H(h(y))=(1-\mu)F(f(x))+\mu F(f(y))\,,\end{array}
$$
that is, the concavity of $F(f)$. Then $f\in \mathcal{C}_\Omega[F]$, and Theorem~\ref{Theorem:1.5} follows.
$\Box$

%%%%%%%%%%%%%%%%%%%
\subsection{The relation between log-concavity and $(1/2)$-log-concavity.}\label{Section:1.5}
%%%%%%%%%%%%%%%%%%%
Next we discuss the relation between $1/2$-log-concavity and log-concavity. 

As we have already said, in our previous paper \cite{IST} our answer to question (Q) was $1/2$-log-concavity.
Indeed, the definition of $F$-concavity given in \cite{IST} is different from here. 
Precisely, let $I=[0,a]$ for some $a>0$ or $I=[0,+\infty)$; following \cite[Definition 2.2]{IST}, given a convex set $\Omega$ and an admissible $F$ on $I$ such that $F(0)=-\infty$, 
we define 
$$
{\mathcal C}_\Omega^*[F]:=
\left\{\,f:\,\mbox{$f$ is a function in $\Omega$ such that 
$\kappa f\in{\mathcal C}_\Omega[F]$ for sufficiently small $\kappa>0$}\,\right\}.
$$
Then:
\begin{itemize}
\item[{\rm (a')}] 
a nonnegative bounded function $f$ is said {\em $F$-concave} in $\Omega$ if $f\in{\mathcal C}_\Omega^*[F]$;
\item[{\rm (b')}] 
let $F_1$ and $F_2$ be admissible on $I_1$ and $I_2$, respectively, where for $i=1,2$ either $I_i=[0,a_i]$ for some $a_i>0$ or $I_i=[0,+\infty)$, then $F_1$-concavity is said  {\em stronger} than $F_2$-concavity in $\Omega$ if ${\mathcal C}_\Omega^*[F_1]\subset{\mathcal C}_\Omega^*[F_2]$.
\end{itemize}
These definitions make $F$-concavity and the comparison between different $F$-concavities independent of the interval $I$.
In this way, thanks to the forthcoming Propositions~\ref{Proposition:4.1} and \ref{Proposition:4.2}, 
we can say that {\em the answer to question {\rm (Q)} is $(1/2)$-log-concavity}.

On the other hand, Definition~{\rm\ref{Definition:1.1}}~{\rm (a)} seems much more natural than 
{\rm (a')} above, and more consistent with the existing literature. Then we preferred to change our notation here, so that question {\rm (Q)} remains strictly dependent on the interval $I$. 
Notice also that in \cite{IST} we used a different terminology for $\alpha$-log-concavity: for $\alpha>0$, a nonnegative bounded function $f$ is said $\alpha$-log-concave in $\Omega$ if $f\in{\mathcal C}_\Omega^*[L_{1/\alpha}]$. 
\medskip

Before giving Proposition \ref{Proposition:4.1}, we characterize $1/2$-log-concavity (compare with Lemma~\ref{Lemma:4.1} and see also Remark \ref{rem-F} (ii)).
\begin{lemma}
\label{Lemma:4.3}
Let $\Omega$ be a convex domain in ${\bf R}^N$. 
Let $a>0$ and $F$ be admissible on $[0,a]$ such that $F$-concavity is nontrivial in $\Omega$. 
Assume that the function
$$
[0,\infty)\ni s\mapsto F\left(a\,e^{-s^2}\right)
$$
is concave. 
Set $\tilde{F}(\tau):=F(a\tau)$ for $\tau \in [0,1]$.
Then $F(0)=-\infty$ and 
$${\mathcal C}_\Omega[L_{1/2}]\subset{\mathcal C}_\Omega[\tilde{F}]=a^{-1}{\mathcal C}_\Omega[F].$$
\end{lemma}
{\bf Proof.}  
Lemma~\ref{Lemma:4.3} follows from the same argument 
as in the proof of \cite[Lemma~2.3]{IST}. 
We give the proof for completeness of this paper. 

Since $F$ is strictly increasing in $(0,1)$, we see that 
the function ${\bf R}\ni s\mapsto F(a\,e^{-s^2})$ is concave and nonconstant, 
which together with Theorem~\ref{Theorem:1.1}~(b) implies that  
$$
\lim_{\tau\to +0}F(\tau)=-\infty=F(0).
$$
Let  $f\in{\mathcal C}_\Omega[L_{1/2}]$. 
Set
$$
g(x):=-\frac12L_{1/2}(f(x))+1
=\left\{
\begin{array}{ll}
\sqrt{-\log f(x)}\quad & \mbox{if}\quad f(x)>0,\vspace{3pt}\\
\infty\quad & \mbox{if}\quad f(x)=0.
\end{array}
\right.
$$
Then $g$ is nonnegative and convex in $\Omega$, that is,
\[
0\le g((1-\mu)x+\mu y)\le(1-\mu)g(x)+\mu g(y)
\]
for $x$, $y\in\Omega$ and $\mu\in[0,1]$.  
On the other hand, 
by $F$-concavity of $ae^{-s^2}$  
we have 
\[
F\left(a\,e^{-\{(1-\mu)g(x)+\mu g(y)\}^2}\right)
\ge(1-\mu)F\left(a\,e^{-g(x)^2}\right)+\mu F\left(a\,e^{-g(y))^2}\right)
\]
for $x$, $y\in\Omega$.   
Combining these with the monotonicity of $F$, 
we obtain 
\begin{equation*}
\begin{split}
 & F(af((1-\mu)x+\mu y))=F\left(a\,e^{-g((1-\mu)x+\mu y)^2}\right)\\
 & \qquad
 \ge F\left(a\,e^{-\{(1-\mu)g(x)+\mu g(y)\}^2}\right)
\ge(1-\mu)F\left(a\,e^{-g(x)^2}\right)+\mu F\left(a\,e^{-g(y)^2}\right)\\
 & \qquad
 =(1-\mu)F(af(x))+\mu F(af(y))
\end{split}
\end{equation*}
if $f(x)f(y)>0$. 
This inequality also holds in the case of $f(x)f(y)=0$. 
Thus $af$ is $F$-concave in $\Omega$ and the proof is complete.
$\Box$

By Lemma~\ref{Lemma:4.3} we apply the argument as in the proof of Lemma~\ref{Lemma:4.1} 
to obtain the following proposition. 
\begin{proposition}
\label{Proposition:4.1}
Let $\Omega$ be a convex domain in ${\bf R}^N$. 
Let $I$ be either $[0,a]$ for some $a>0$ or $[0,+\infty)$, and let $F$ be admissible on $I$ such that $F(0)=-\infty$. 
Assume that 
\begin{equation}
\label{eq:4.11}
e^{t\Delta_\Omega}\varphi\in{\mathcal C}_\Omega^*[F]\quad
\mbox{for $\varphi\in {\mathcal C}_\Omega^*[F]\cap L^\infty(\Omega)$ and $t>0$}.
\end{equation}
Then ${\mathcal C}_\Omega^*[L_{1/2}]\subset{\mathcal C}^*_\Omega[F]$. 
\end{proposition}
{\bf Proof.} 
Proposition~\ref{Proposition:4.1} is a slight modification of \cite[Theorem~3.2]{IST}. 
By the same argument as in the proof of Lemma~\ref{Lemma:4.1} 
we see that the function
$$
[0,\infty)\ni s\mapsto F\left(ke^{-s^2}\right)
$$
is concave for sufficiently small $k>0$. 
Set $F_k(\tau):=F(k\tau)$ for $\tau\in[0,1]$. 
Then we deduce from Lemma~\ref{Lemma:4.3} that ${\mathcal C}_\Omega[L_{1/2}]\subset{\mathcal C}_\Omega[F_k]=k^{-1}\mathcal{C}_\Omega[F]$,
whence ${\mathcal C}_\Omega^*[L_{1/2}]\subset{\mathcal C}^*_\Omega[F]$, as desired.~$\Box$
\vspace{3pt}

Furthermore, by \cite[Lemma~2.1, Theorem~3.1]{IST} we have:
\begin{proposition}
\label{Proposition:4.2}
Let $\Omega$ be a convex domain in ${\bf R}^N$. 
Then, for any $1/2\le\alpha\le 1$, 
$\alpha$-log-concavity is preserved by the Dirichlet heat flow in $\Omega$. Furthermore, 
$$
\bigcup_{t>0}e^{t\Delta_\Omega}{\mathcal C}_\Omega^*[L_\alpha]\subset {\mathcal C}_\Omega^*[L_\alpha].
$$
\end{proposition}
We observe from Propositions~\ref{Proposition:4.1} and \ref{Proposition:4.2} 
that $1/2$-log-concavity is the strongest 
among $F$-concavities satisfying $F(0)=-\infty$ and 
$$
\bigcup_{t>0}e^{t\Delta_\Omega}{\mathcal C}_\Omega^*[F]\subset {\mathcal C}_\Omega^*[F].
$$

On the other hand, set
$$
L_{1/2}^k(\tau):=L_{1/2}(k^{-1}\tau)\quad\mbox{for}\quad\tau\in[0,k],
$$
where $k>0$. 
Then $L_{1/2}^k$ is admissible on $[0,k]$. 
By Example~\ref{Example:1.3}~(iii) we see that 
$$
{\mathcal C}_\Omega[L_{1/2}^{k_1}]\subset {\mathcal C}_\Omega[L_{1/2}^{k_2}]
\quad\mbox{if}\quad k_1\le k_2.
$$
For $k>1$ we normalize $L_{1/2}^k$ as follows:
$$
{\mathcal L}_{1/2}^k(\tau):=(\log k)^{\frac{1}{2}}\left( L_{1/2}^k(\tau)-L_{1/2}^k(1)\right)
\quad \mbox{for}\quad\tau\in[0,k].
$$
Then $({\mathcal L}_{1/2}^k)(1)=0$ and $({\mathcal L}_{1/2}^k)'(1)=1$. 
Furthermore, ${\mathcal C}_\Omega[L^k_{1/2}]={\mathcal C}_\Omega[{\mathcal L}^k_{1/2}]$ 
and 
\begin{equation*}
\begin{split}
{\mathcal L}_{1/2}^k(\tau)
 & =-2(\log k)^{\frac{1}{2}}\left\{ (-\log\tau+\log k)^{\frac{1}{2}}-(\log k)^{\frac{1}{2}}\right\}\\
 & =-2(\log k)\left\{ (1-(\log k)^{-1}\log\tau)^{\frac{1}{2}}-1\right\}
 \to \log\tau
\end{split}
\end{equation*}
as $k\to\infty$ for $\tau>0$. These mean that 
log-concavity can be regarded as the limit of $L_{1/2}^k$-concavity despite 
$$
\bigcup_{k>0}{\mathcal C}_\Omega[L_{1/2}^k]\subsetneq{\mathcal C}_\Omega[L_0].
$$

%%%%%%%%%%%%%%%%%%%%%%%%%%%%%%%%%%%%
%%%%%%%%%%%%%%%%%%%%%%%%%%%%%%%%%%%%
\section{Open problems}
%%%%%%%%%%%%%%%%%%%%%%%%%%%%%%%%%%%%
%%%%%%%%%%%%%%%%%%%%%%%%%%%%%%%%%%%%
We present some open problems related to the Dirichlet heat flow.
Let $\Omega$ be a bounded convex domain in ${\bf R}^N$. 
Let $\varphi$ be a (nontrivial) continuous function in $\Omega$ 
such that $0\le\varphi\le 1$ in $\Omega$. 
Applying the standard theory for parabolic equations, 
we see that 
\begin{equation}
\label{eq:5.1}
\lim_{t\to\infty}e^{\lambda_1 t}e^{t\Delta_\Omega}\varphi=(\varphi,\phi)_{L^2(\Omega)}\phi
\end{equation}
in the sense of $C^2(\overline{\Omega})$. 
Here $\phi$ is the first normalized Dirichlet eigenfunction in $\Omega$ and $\lambda_1$ 
the first Dirichlet eigenvalue, that is, $\phi$ satisfies $\|\phi\|_{L^2(\Omega)}=1$ and 
$$
\left\{
\begin{array}{ll}
-\Delta\phi=\lambda_1\phi & \mbox{in}\quad\Omega\,,\\
\phi>0 & \mbox{in}\quad\Omega\,,\\
\phi=0 & \mbox{on}\quad\partial\Omega\,.
\end{array}
\right.
$$
Assume that $\varphi$ is $1/2$-log-concave in $\Omega$. 
Although it follows from Proposition~\ref{Proposition:4.1} that 
$e^{t\Delta_\Omega}\varphi$ is $1/2$-log-concave in $\Omega$ for $t>0$, 
the limit procedure in \eqref{eq:5.1} does not imply 
that $\phi$ is $1/2$-log-concave in $\Omega$ even if $\phi\le 1$ in $\Omega$. 
Indeed, $1/2$-log-concavity is not closed under positive scalar multiplication (see Theorem~\ref{Theorem:1.3}). 
Despite this, we conjecture 
that the first positive Dirichlet eigenfunction $\phi$ of the Laplacian in $\Omega$
is $1/2$-log-concave if $\phi\le 1$ in $\Omega$.

In the same spirit, it would be interesting to determine the weakest $F$-concavity 
preserved by the Dirichlet heat flow. 
In \cite[Theorem~4.1]{IS2} it is showed that 
$p$-concavity is not preserved by the Dirichlet heat flow in $\Omega$ 
for some $p\in(-\infty,0)$. More precisely, 
for any $t_*>0$,  
there exists $\varphi\in C_0(\Omega)$ with the following properties:
\begin{itemize}
  \item $\varphi$ is $p$-concave in $\Omega$ for some $p\in(-\infty,0)$;
  \item $u(t)=e^{t\Delta_\Omega}\varphi$ is not quasiconcave in $\Omega$, in particular, not $p$-concave, 
  for some $t\in(0,t_*)$.
\end{itemize}
See also \cite{IS1}. 
However, the identification of the weakest $F$-concavity (or even just about the weakest power concavity) 
preserved by the Dirichlet heat flow is still open. 
\medskip

\noindent
{\bf Acknowledgements.} 
The first author was supported in part by the Grant-in-Aid for Scientific Research (S)(No.~19H05599)
from Japan Society for the Promotion of Science.
The second author has been partially supported by INdAM through a GNAMPA Project.
The third author was supported in part by the Grant-in-Aid for Scientific Research (C) (No.~19K03494).
%%%%%%%%%%%%%%%%%%%%%%%%%%%%%%%%%%%%%%
%%%%%%%%%%%%    references    %%%%%%%%%%%%%%%%%%
%%%%%%%%%%%%%%%%%%%%%%%%%%%%%%%%%%%%%%

%%%%%%%%%%%%%%%%%%%%%%%%%%%%%%%%%%%%
%%%%%%%%%%%%%%%%%%%%%%%%%%%%%%%%%%%%

\begin{bibdiv}
  \begin{biblist}
\bib{An}{article}{
   author={An, Mark Yuying},
   title={Logconcavity versus logconvexity: a complete characterization},
   journal={J. Econom. Theory},
   volume={80},
   date={1998},
   pages={350--369},
}
%%%%%%%
\bib{BL}{article}{
   author={Brascamp, Herm Jan},
   author={Lieb, Elliott H.},
   title={On extensions of the Brunn-Minkowski and Pr\'{e}kopa-Leindler
   theorems, including inequalities for log concave functions, and with an
   application to the diffusion equation},
   journal={J. Functional Analysis},
   volume={22},
   date={1976},
   pages={366--389},
}
%%%%%%
\bib{Co}{article}{
   author={Colesanti, Andrea},
   title={Log-concave functions},
   conference={
      title={Convexity and concentration},
   },
   book={
      series={IMA Vol. Math. Appl.},
      volume={161},
      publisher={Springer, New York},
   },
   date={2017},
   pages={487--524},
}
%%%%%%
%%%%%%
\bib{G}{article}{
   author={Gabriel, R. M.},
   title={A result concerning convex level surfaces of $3$-dimensional
   harmonic functions},
   journal={J. London Math. Soc.},
   volume={32},
   date={1957},
   pages={286--294},
}
%%%%%%
\bib{HLP}{book}{
   author={Hardy, G. H.},
   author={Littlewood, J. E.},
   author={P\'{o}lya, G.},
   title={Inequalities},
   series={Cambridge Mathematical Library},
   note={Reprint of the 1952 edition},
   publisher={Cambridge University Press, Cambridge},
   date={1988},
   pages={xii+324},
}
%%%%%%
\bib{IS1}{article}{
   author={Ishige, Kazuhiro},
   author={Salani, Paolo},
   title={Is quasi-concavity preserved by heat flow?},
   journal={Arch. Math. (Basel)},
   volume={90},
   date={2008},
   pages={450--460},
 }
%%%%%
\bib{IS2}{article}{
   author={Ishige, Kazuhiro},
   author={Salani, Paolo},
   title={Convexity breaking of the free boundary for porous medium
   equations},
   journal={Interfaces Free Bound.},
   volume={12},
   date={2010},
   pages={75--84},
}
%%%%%
\bib{IST}{article}{
   author={Ishige, Kazuhiro},
   author={Salani, Paolo},
   author={Takatsu, Asuka},
   title={To logconcavity and beyond},
   journal={Commun. Contemp. Math.},
   volume={22},
   date={2020},
   number={2},
   pages={1950009, 17},
}
\bib{ML}{article}{
   author={Makar-Limanov, L. G.},
   title={The solution of the Dirichlet problem for the equation $\Delta
   u=-1$ in a convex region},
   journal={Mat. Zametki},
   volume={9},
   date={1971},
   pages={89--92},
}
\bib{Na}{book}{
   author={Naudts, Jan},
   title={Generalised thermostatistics},
   publisher={Springer-Verlag London, Ltd., London},
   date={2011},
   pages={x+201},
}
%%%%%%%%
\bib{SW}{article}{
   author={Saumard, Adrien},
   author={Wellner, Jon A.},
   title={Log-concavity and strong log-concavity: a review},
   journal={Stat. Surv.},
   volume={8},
   date={2014},
   pages={45--114},
}
%%%%%%%%
\bib{Sh}{book}{
   author={Schneider, Rolf},
   title={Convex bodies: the Brunn-Minkowski theory},
   series={Encyclopedia of Mathematics and its Applications},
   volume={44},
   publisher={Cambridge University Press, Cambridge},
   date={1993},
   pages={xiv+490},
   isbn={0-521-35220-7},
   doi={10.1017/CBO9780511526282},
}
%%%%%%%%
\end{biblist}
\end{bibdiv}
\end{document}